\documentclass[final]{siamltex}
\usepackage{amsmath}
\usepackage{amsfonts}
\usepackage{graphicx}
\usepackage{subfigure}
\def\fig#1{Figure~\ref{#1}}

\def\tab#1{Table~\ref{#1}}

\def\bE{\mathbf E}
\def\bH{\mathbf H}
\def\bJ{\mathbf J}

\def\bx{\mathbf x}

\def\bk{\mathbf k}

\def\b0{\mathbf 0}

\title{Transverse electric scattering on inhomogeneous objects: 
spectrum of integral operator and preconditioning}

\author{Grigorios~P.~Zouros\thanks{School of Electrical and Computer 
        Engineering, National Technical University of Athens, 
        Athens 15773, Greece ({\tt zouros@mail.ntua.gr}).}
        \and Neil~V.~Budko\thanks{Numerical Analysis, DIAM, 
        Faculty of Electrical Engineering, Mathematics, 
        and Computer Science, Delft University of Technology, 2628 CD Delft, 
        The Netherlands ({\tt n.v.budko@tudelft.nl}).}}

\begin{document}

\maketitle

\begin{abstract}
The domain integral equation method with its FFT-based matrix-vector products is a viable alternative 
to local methods in free-space scattering problems. However, it often suffers from the extremely 
slow convergence of iterative methods, especially 
in the transverse electric (TE) case with large or negative permittivity.
We identify the nontrivial essential spectrum of the pertaining integral 
operator as partly responsible for this behavior, and the main  
reason why a normally efficient deflating preconditioner does not work. We solve this 
problem by applying an explicit multiplicative regularizing operator, which transforms 
the system to the form `identity plus compact', yet allows the resulting matrix-vector products
to be carried out at the FFT speed. Such a regularized system is then further preconditioned
by deflating an apparently stable set of eigenvalues with largest magnitudes, which results in a robust
acceleration of the restarted GMRES under constraint memory conditions. 
\end{abstract}

\begin{keywords} 
Domain integral equation, singular integral operators, electromagnetism, 
TE scattering, spectrum of operators, essential spectrum, regularizer, 
deflation, preconditioner
\end{keywords}

\begin{AMS}
78A45, 65F08, 45E10, 47G10, 15A23
\end{AMS}

\pagestyle{myheadings}
\thispagestyle{plain}
\markboth{G.~P.~ZOUROS AND N.~V.~BUDKO}{DOMAIN INTEGRAL OPERATOR: SPECTRUM AND PRECONDITIONING}

\section{Introduction}
There is a steady interest in the numerical simulation of the electromagnetic field in inhomogeneous media.
The methods can be roughly divided into two categories -- local and global -- in accordance with the governing
equations. The local methods based on the differential Maxwell's equations 
\cite{Erlangga_2008, Erlangga_2006, ihlenburg_Babuska_1995a, ihlenburg_Babuska_1995b} are 
generally more popular due to the sparse nature of their matrices and the ease of programming. 
Alongside the local methods, global methods based on an equivalent integral 
equation formulation are also often employed, especially in free-space scattering.
If an object is large and homogeneous or has a perfectly conducting boundary then the problem is usually reduced
to a boundary integral equation with the fields (currents) at the interfaces being the fundamental 
unknowns \cite{Harrington_1993}. If an object is continuously inhomogeneous or is a composite 
consisting of many small different parts, the most appropriate global method is the domain integral equation (DIE),
in two dimensions, or the volume integral equation, in three dimensions 
\cite{Ilinski_2000, Joachimowicz_Pichot_1990, Kleinman_Roach_VanDenBerg_1990, Livesay_Chen_1974, 
Peterson_Klock_1988, Richmond_1966, Sondergaard_2007,  Su_1987, Volakis_1992}. 
Although the global methods produce dense matrices, they are generally more stable with respect to 
discretization than the local ones, and the convolution-type integral operators 
sometimes allow to compute matrix-vector products at the FFT speed. 
These properties make the DIE method a viable alternative to local methods for 
certain free-space scattering problems.

The main difficulties with the DIE method are the non-normality of both the operator and 
the resulting system matrix, inherent to frequency-domain electromagnetic scattering, 
and the extremely slow convergence of the few iterative methods that
can be applied with such matrices. Numerical experiments consistently show that the GMRES
algorithm is the best choice  \cite{Fan_2006}, and it can be proven that the full (un-restarted) 
GMRES will evetually converge 
with any physically meaningfull scatterer and incident field \cite{SamokhinBook}.
Most of the practically interesting problems, however, involve the number of unknowns prohibitive for
the full GMRES, and its restarted version is used. 
Unfortunately, sometimes, especially in TE scattering with large permittivities/dimensions, 
restarted GMRES as well Bi-CGSTAB converge very slowly or do not converge at all.

Several successful preconditioning strategies have been proposed for local methods, 
see e.g. \cite{Erlangga_2008, Erlangga_2006}, and boundary integral formulations
\cite{Andriulli_2008, Bagci_2010, Christiansen_2003, Coifman_1993, Song_1997}. 
Whereas the domain integral equation method is
lagging behind in this respect \cite{Card_2004, Fan_2006, Flatau_1997, Yurkin_2007},
reporting convergence and acceleration thereof mainly for low-contrast and/or small objects.
One of the difficulties in designing a suitable multiplicative preconditioner for this method is 
that it needs to be either sparse or have a block-Toeplitz form to be able to compete with 
local methods in terms of memory and speed.
Recently a deflation-based preconditioner has been proposed in \cite{sifuentes} for accelerating the
DIE solver in the TM case. A thorough spectral analysis of the system matrix reported in 
\cite{sifuentes} showed that the outlying eigenvalues (largest in magnitude) are responsible for a 
period of stagnation of the full GMRES, and that their deflation accelerates the iterative process.
As we show here this strategy, however, does not work in the TE case. Moreover, using the  
deflation on top of a full GMRES puts an even higher strain on the memory.

To further extend the applicability of the DIE 
method towards objects with higher contrasts/sizes, in this paper we 
propose a two-stage preconditioner based on the regularization of the pertaining singular integral operator
and subsequent deflation of the largest eigenvalues. To derive a regularizer we employ 
the symbol calculus, which also provides us with the essential spectrum of the original operator.
Our regularizer may be viewed as a generalization of the Calder{\'o}n identity recently employed for
preconditioning of the boundary integral equation method \cite{Andriulli_2008, Bagci_2010, Christiansen_2003}.
As an illustration and a proof of the concept we compute the complete spectrum of the system matrix
for a few typical objects of resonant size. We
analyze the difference in spectra between the TE and TM polarizations and demonstrate
that the discretized version of the regularizer indeed contracts a very dense group of eigenvalues making
the regularized system amenable to deflation. Unlike the previously mentioned boundary integral
equation approach \cite{Andriulli_2008} 
the DIE method permits a straightforward discretization of the
continuous regularizer without any additional precautions. As our regularizer is a convolution-type
integral operator, all matrix-vector products can be carried out with the FFT algorithm.
The subsequent deflation of the largest eigenvalues can be achieved
with the standard Matlab's
{\tt eigs} routine. Most importantly, these largest eigenvalues and the assiociated eigenvectors 
turn out to be rather stable, so that one can considerably accelerate the {\tt eigs} algorithm by 
setting a modest tolerance. We have performed a number of
numerical experiments which have consistently demonstrated that given a fixed amount of computer
memory it is better to spend some part of it on deflation than all of it on maximizing the dimension of 
the inner Krylov subspace of the restarted GMRES. Keeping the {\tt restart} parameter of the GMRES
in the order or larger than the dimension of the deflation subspace gives a speed-up which grows with the object 
size and permittivity.

Although the DIE method (both TM and TE versions of it) has been around in engineering 
community for many years 
\cite{Ilinski_2000, Joachimowicz_Pichot_1990, Kleinman_Roach_VanDenBerg_1990, Livesay_Chen_1974, 
Peterson_Klock_1988, Richmond_1966, Sondergaard_2007,  Su_1987, Volakis_1992}, 
the full theoretical analysis of the two-dimensional case (in the strongly singular form)
is still missing. The reason might be the presence of special functions in the kernel of the integral operator,
which makes calculations more tedious than in the three-dimensional case \cite{bud_sam_06, SamokhinBook}. 
Therefore, before going into the numerical details, we devote some space here to the analysis of the 
two-dimensional singular integral operator. As a result we show an almost complete equivalence of 
the TE scattering to the more general three-dimensional scattering problem. 
Explicit numerical calculations of the matrix spectrum for resonant scatterers presented here 
are generally not possible with three-dimensional objects of that size. Whereas, the established 
spectral equivalence indicates that the proposed preconditioner will work with 
three-dimensional problems as well. 

\section{Domain integral equation}
We consider the two-dimensional frequency-domain electromagnetic 
scattering problem on an inhomogeneous object of finite spatial 
support exhibiting contrast in both electric and magnetic properties. 
We assume that the object (a cylinder) 
as well the source of the incident field (a plane wave, current carrying line, etc.) 
are invariant in the $x_{3}$-direction of the spatial Cartesian system.
The Maxwell equations describing the total electromagnetic field in the 
$(x_{1},x_{2})$-plane in the presence of an inhomogeneous
magneto-dielectric scatterer can be written as
\begin{align}
\label{eq:MaxwellTE}
\begin{split}
-\partial_{2}H_{3}(\bx,\omega)-i\omega\varepsilon(\bx,\omega) E_{1}(\bx,\omega) & = -J_{1}(\bx,\omega),
\\
\partial_{1}H_{3}(\bx,\omega)-i\omega\varepsilon(\bx,\omega) E_{2}(\bx,\omega) & = -J_{2}(\bx,\omega),
\\
-\partial_{2}E_{1}(\bx,\omega)+\partial_{1}E_{2}(\bx,\omega)-i\omega\mu(\bx,\omega) H_{3}(\bx,\omega) & = 0,
\end{split}
\\
\label{eq:MaxwellTM}
\begin{split}
\partial_{2}H_{1}(\bx,\omega)-\partial_{1}H_{2}(\bx,\omega)-i\omega\varepsilon(\bx,\omega) E_{3}(\bx,\omega) & = -J_{3}(\bx,\omega),
\\
\partial_{2}E_{3}(\bx,\omega)-i\omega\mu(\bx,\omega) H_{1}(\bx,\omega) & = 0,
\\
-\partial_{1}E_{3}(\bx,\omega)-i\omega\mu(\bx,\omega) H_{2}(\bx,\omega) & = 0,
\end{split}
\end{align}
where we introduce the two-dimensional position vector $\bx=(x_{1},x_{2})$;
$\bJ$ is the electric current density; 
$\varepsilon$ 
and $\mu$ are the possibly complex-valued
dielectric permittivity and magnetic permeability; and $\omega$
is the angular frequency. It is easy to notice that the first three equations, (\ref
{eq:MaxwellTE}), are decoupled from the last three, (\ref{eq:MaxwellTM}). Reflecting the fact that 
the electric field strength is {\it transverse} to the axis of invariance, 
equations (\ref{eq:MaxwellTE}) are said to describe 
the transverse electric or TE case (also somewhat confusingly called H-wave polarization). 
Analogously, the three equations (\ref{eq:MaxwellTM})
correspond to the transverse magnetic, or TM case (also known as E-wave polarization).

The integral formulation of the scattering problem is obtained by first introducing
the so-called incident field $(\bE^{\rm in}$, $\bH^{\rm in})$, which has the same physical 
source as the total field, but satisfies some simpler version of the above equations. 
Usually, the incident field is considered in a homogeneous 
isotropic background medium with parameters $\varepsilon_{\rm b}(\omega)$ 
and $\mu_{\rm b}(\omega)$. Then, it is fairly easy to see that the scattered field,
defined as the difference between the total and the incident fields
\begin{align}
\label{eq:ScatteredField}
\bE^{\rm sc}(\bx,\omega)=\bE(\bx,\omega)-\bE^{\rm in}(\bx,\omega),
\;\;\;
\bH^{\rm sc}(\bx,\omega)=\bH(\bx,\omega)-\bH^{\rm in}(\bx,\omega),
\end{align}
also satisfies Maxwell's equations in the same homogeneous medium, i.e.,
\begin{align}
\label{eq:MaxwellTEsc}
\begin{split}
-\partial_{2}H^{\rm sc}_{3}(\bx)-i\omega\varepsilon_{\rm b} E^{\rm sc}_{1}(\bx) & = -J^{\rm ind}_{1}(\bx),
\\
\partial_{1}H^{\rm sc}_{3}(\bx)-i\omega\varepsilon_{\rm b} E^{\rm sc}_{2}(\bx) & = -J^{\rm ind}_{2}(\bx),
\\
-\partial_{2}E^{\rm sc}_{1}(\bx)+\partial_{1}E^{\rm sc}_{2}(\bx)-i\omega\mu_{\rm b} H^{\rm sc}_{3}(\bx) & = -K^{\rm ind}_{3}(\bx),
\end{split}
\\
\label{eq:MaxwellTMsc}
\begin{split}
\partial_{2}H^{\rm sc}_{1}(\bx)-\partial_{1}H^{\rm sc}_{2}(\bx)-i\omega\varepsilon_{\rm b} E^{\rm sc}_{3}(\bx) & = -J^{\rm ind}_{3}(\bx),
\\
\partial_{2}E^{\rm sc}_{3}(\bx)-i\omega\mu_{\rm b} H^{\rm sc}_{1}(\bx) & = -K^{\rm ind}_{1}(\bx),
\\
-\partial_{1}E^{\rm sc}_{3}(\bx)-i\omega\mu_{\rm b} H^{\rm sc}_{2}(\bx) & = -K^{\rm ind}_{2}(\bx),
\end{split}
\end{align}
where the implicit dependence on $\omega$ has been omitted for simplicity. The physical 
sources of the scattered field are the {\it induced} current densities, also known as contrast currents:
\begin{align}
\label{eq:ContrastCurrents}
\begin{split}
J^{\rm ind}_{k}(\bx)&=-i\omega\left[\varepsilon(\bx)-\varepsilon_{\rm b}\right]E_{k}(\bx), \;\;\;k=1,2,3;
\\
K^{\rm ind}_{m}(\bx)&=-i\omega\left[\mu(\bx)-\mu_{\rm b}\right]H_{m}(\bx), \;\;\;m=1,2,3.
\end{split}
\end{align}
It is convenient to write the Maxwell equations in the following matrix form:
\begin{align}
\label{eq:MatrixTE}
\begin{bmatrix}
-i\omega\varepsilon_{\rm b} & 0 & -\partial_{2}\\
0 & -i\omega\varepsilon_{\rm b} & \partial_{1} \\
-\partial_{2} & \partial_{1} & -i\omega\mu_{\rm b}
\end{bmatrix}
\begin{bmatrix}
E_{1}^{\rm sc}\\
E_{2}^{\rm sc}\\
H_{3}^{\rm sc}
\end{bmatrix}
&=
\begin{bmatrix}
-J_{1}^{\rm ind}\\
-J_{2}^{\rm ind}\\
-K_{3}^{\rm ind}
\end{bmatrix},
\\
\label{eq:MatrixTM}
\begin{bmatrix}
-i\omega\mu_{\rm b} & 0 & \partial_{2}\\
0 & -i\omega\mu_{\rm b} & -\partial_{1} \\
\partial_{2} & -\partial_{1} & -i\omega\varepsilon_{\rm b}
\end{bmatrix}
\begin{bmatrix}
H_{1}^{\rm sc}\\
H_{2}^{\rm sc}\\
E_{3}^{\rm sc}
\end{bmatrix}
&=
\begin{bmatrix}
-K_{1}^{\rm ind}\\
-K_{2}^{\rm ind}\\
-J_{3}^{\rm ind}
\end{bmatrix}.
\end{align}
The solution of these equations can be obtained by first transforming 
(\ref{eq:MatrixTE})--(\ref{eq:MatrixTM}) to the $\bk$-domain with the help of the 
two-dimensional spatial Fourier transform, then solving the resulting algebraic 
systems, and finally transforming the result back to the $\bx$-domain. 
The obtained solution expresses the scattered field in terms of the induced currents
and, via (\ref{eq:ContrastCurrents}), in terms of the unknown total field. Subsequently, using 
(\ref{eq:ScatteredField}), one can replace the scattered field with the difference between 
the total and the incident fields, thus arriving at the following integro-differential equations 
with the total field as a fundamental unknown:
\begin{align}
\label{eq:MatrixSolutionTE}
\begin{bmatrix}
E_{1}^{\rm in}\\
E_{2}^{\rm in}\\
H_{3}^{\rm in}
\end{bmatrix}
&=
\begin{bmatrix}
E_{1}\\
E_{2}\\
H_{3}
\end{bmatrix}
-
\begin{bmatrix}
k_{\rm b}^{2}+\partial_{1}^{2} & \partial_{1}\partial_{2} & 
-i\omega\mu_{\rm b}(-\partial_{2})\\
\partial_{2}\partial_{1} & k_{\rm b}^{2}+\partial_{2}^{2} &
-i\omega\mu_{\rm b}\partial_{1}\\
-i\omega\varepsilon_{\rm b}(-\partial_{2}) & -i\omega\varepsilon_{\rm b}\partial_{1} &
k_{\rm b}^{2}
\end{bmatrix}
\begin{bmatrix}
g*(\chi_{{\rm e}}E_{1})\\
g*(\chi_{{\rm e}}E_{2})\\
g*(\chi_{{\rm m}}H_{3})
\end{bmatrix},
\\
\label{eq:MatrixSolutionTM}
\begin{bmatrix}
H_{1}^{\rm in}\\
H_{2}^{\rm in}\\
E_{3}^{\rm in}
\end{bmatrix}
&=
\begin{bmatrix}
H_{1}\\
H_{2}\\
E_{3}
\end{bmatrix}
-
\begin{bmatrix}
k_{\rm b}^{2}+\partial_{1}^{2} & \partial_{1}\partial_{2} & 
-i\omega\varepsilon_{\rm b}\partial_{2}\\
\partial_{2}\partial_{1} & k_{\rm b}^{2}+\partial_{2}^{2} &
-i\omega\varepsilon_{\rm b}(-\partial_{1})\\
-i\omega\mu_{\rm b}\partial_{2} & -i\omega\mu_{\rm b}(-\partial_{1}) &
k_{\rm b}^{2}
\end{bmatrix}
\begin{bmatrix}
g*(\chi_{{\rm m}}H_{1})\\
g*(\chi_{{\rm m}}H_{2})\\
g*(\chi_{{\rm e}}E_{3})
\end{bmatrix}.
\end{align}
Here we have introduced the wavenumber of the homogeneous background medium
$k_{\rm b}^{2}=\omega^{2}\varepsilon_{\rm b}\mu_{\rm b}$ and the normalized 
electric and magnetic contrast functions:
\begin{align}
\label{eq:ElectricContrast}
\chi_{\rm {\rm e}}(\bx,\omega)&=\frac{\varepsilon(\bx,\omega)}{\varepsilon_{\rm b}}-1,
\\
\chi_{\rm {\rm m}}(\bx,\omega)&=\frac{\mu(\bx,\omega)}{\mu_{\rm b}}-1.
\end{align}
The star $(*)$ denotes the two-dimensional convolution, e.g.,
\begin{align}
\label{eq:Convolution}
g*(\chi_{{\rm e}}E_{1})= \int_{\bx'\in{\mathbb R}^{2}}\,
g(\bx-\bx')\chi_{{\rm e}}(\bx')E_{1}(\bx'){\rm d}\bx',
\end{align}
The scalar Green's function of the two-dimensional Helmholtz equation
is given by
\begin{align}
\label{eq:GreenFunction}
g(\bx,\omega)=\frac{i}{4}H^{(1)}_{0}(k_{\rm b}\vert\bx\vert),
\end{align}
where $H^{(1)}_{0}$ is the zero's order Hankel function of the first kind.

\section{Operator, symbol, spectrum, and regularizer}\label{op_sym_sp_reg}
Although the integro-differential equations
(\ref{eq:MatrixSolutionTE})--(\ref{eq:MatrixSolutionTM}) can to a certain extent 
be analyzed directly on the pertaining Sobolev spaces, we prefer to transform them into
singular integral equations and work on the Hilbert space $L^{(3)}_{2}(D)$ of 
vector-valued functions with spatial support on $D$.

As the two problems (\ref{eq:MatrixSolutionTE}) and (\ref{eq:MatrixSolutionTM}) 
are mathematically identical, up to some constants and non-essential sign changes in 
the differential matrix operator, we shall concentrate on just one of them,
say the TE case (\ref{eq:MatrixSolutionTE}). A standard singular integral equation
is obtained from (\ref{eq:MatrixSolutionTE}) by carrying out the spatial 
derivatives of the weakly singular integrals (\ref{eq:Convolution}),
and writing the result as
\begin{align}
\label{eq:StandardForm}
u^{\rm in}=Du+GXu+KXu,
\end{align}
where $D$ and $X$ are ``diagonal'' multiplication operators, $G$ is
a principal-value singular integral operator, and $K$ is a compact
integral operator. In our case, things become a little more complicated
due to the matrix-valued kernel and the involvement of special functions,
but eventually one can put equation (\ref{eq:MatrixSolutionTE}) in the standard form as
\begin{align}
\label{eq:SingularTE}
\begin{split}
\begin{bmatrix}
E_{1}^{\rm in}\\
E_{2}^{\rm in}\\
H_{3}^{\rm in}
\end{bmatrix}
=
\begin{bmatrix}
S & 0 & 0 \\
0 & S & 0 \\
0 & 0 & I
\end{bmatrix}
\begin{bmatrix}
E_{1}\\
E_{2}\\
H_{3}
\end{bmatrix}
&
+
p.\,v.
\begin{bmatrix}
G_{11} & G_{12} & 0 \\
G_{21} & G_{22} & 0 \\
0 & 0 & 0
\end{bmatrix}*
\begin{bmatrix}
X_{{\rm e}}E_{1}\\
X_{{\rm e}}E_{2}\\
X_{{\rm m}}H_{3}
\end{bmatrix}
\\
&+
\begin{bmatrix}
K_{11} & K_{12} & K_{13} \\
K_{21} & K_{22} & K_{23} \\
K_{31} & K_{32} & K_{33}
\end{bmatrix}*
\begin{bmatrix}
X_{{\rm e}}E_{1}\\
X_{{\rm e}}E_{2}\\
X_{{\rm m}}H_{3}
\end{bmatrix}.
\end{split}
\end{align}
Here, the first matrix operator on the right contains the identity operator $I$,
and two operators $S$, which denote pointwise multiplication with the following function:
\begin{align}
\label{eq:FreeTerms}
s(\bx)=1+\frac{1}{2}\chi_{{\rm e}}(\bx).
\end{align}
Operators $X_{\rm e/m}$ denote pointwise multiplication with
the corresponding contrast functions $\chi_{\rm e/m}(\bx)$.
The second term on the right is the principal value of a
convolution-type integral operator (with circular exclusion area around the singularity, 
which tends to zero in the limit), where the nonzero elements of the matrix-valued kernel are
\begin{align}
\label{eq:PVKernel}
G_{nm}(\bx)=-\frac{1}{2\pi\vert\bx\vert^{2}}
\left[2\theta_{n}\theta_{m}-\delta_{nm}\right],\;\;\;n,m=1,2;
\end{align}
with $\theta_{n}=x_{n}/\vert\bx\vert$ being a Cartesian component of a 
unit vector. The last term on the right in (\ref{eq:SingularTE})
is also a convolution with the matrix elements of the kernel
given by
\begin{align}
\label{eq:KompactKernel1122}
\begin{split}
K_{nm}=&\left[\frac{1}{2\pi\vert\bx\vert^2}-\frac{ik_{\rm b}}{4\vert\bx\vert}H_{1}^{(1)}(k_{\rm b}\vert\bx\vert)
\right]
\left[2\theta_{n}\theta_{m}-\delta_{nm}\right]
\\
&+\frac{ik_{\rm b}^2}{4}H_{0}^{(1)}(k_{\rm b}\vert\bx\vert)
\left[\theta_{n}\theta_{m}-\delta_{nm}\right], \;\;\;n,m=1,2;
\end{split}
\\
\label{eq:KompactKernel13}
\begin{split}
K_{13}=-\frac{\omega\mu_{\rm b}k_{\rm b}\theta_{2}}{4}H_{1}^{(1)}(k_{\rm b}\vert\bx\vert),
\end{split}
\\
\label{eq:KompactKernel31}
\begin{split}
K_{31}=-\frac{\omega\varepsilon_{\rm b}k_{\rm b}\theta_{2}}{4}H_{1}^{(1)}(k_{\rm b}\vert\bx\vert),
\end{split}
\\
\label{eq:KompactKernel23}
\begin{split}
K_{23}=\frac{\omega\mu_{\rm b}k_{\rm b}\theta_{1}}{4}H_{1}^{(1)}(k_{\rm b}\vert\bx\vert),
\end{split}
\\
\label{eq:KompactKernel32}
\begin{split}
K_{32}=\frac{\omega\varepsilon_{\rm b}k_{\rm b}\theta_{1}}{4}H_{1}^{(1)}(k_{\rm b}\vert\bx\vert),
\end{split}
\\
\label{eq:KompactKernel33}
\begin{split}
K_{33}=-\frac{ik_{\rm b}^{2}}{4}H_{0}^{(1)}(k_{\rm b}\vert\bx\vert).
\end{split}
\end{align}

As can be deduced from asymptotic expansions of special 
functions, all $K_{nm}$ kernels are weakly singular, i.e.,
their singularity is less than that of the factor $\vert\bx-\bx'\vert^{-2}$,
and the corresponding integrals exist in the usual sense.
It is also clear that the matrix integral operator $KX$ is compact on $L^{(3)}_{2}(D)$.
The remaining strong singularity is explicitly contained within the
$G_{nm}$ kernel. The latter kernel also features a very special
tensor $2\theta_{n}\theta_{m}-\delta_{nm}$, which guarantees the 
existence of the integral in the sense
of principal value. The standard assumption, sufficient to guarantee the
boundedness of the matrix singular integral operator $G$ and the equivalence of 
the integral equation (\ref{eq:StandardForm}) to the original 
Maxwell's equations, is that of H{\"o}lder continuity  
of the dielectric permittivity and of all components of the incident field vector 
for all $\bx\in{\mathbb R}^{2}$. Although in practice the equation
is posed and solved on $L^{(3)}_{2}(D)$ only, the existence theory was
developed for domains without an edge \cite{Mikhlin}, i.e., 
on $L^{(3)}_{2}({\mathbb R}^{2})$. Samokhin, however, has shown, \cite{SamokhinBook},
that as far as the question of existence is concerned all the operators 
can be extended to $L^{(3)}_{2}({\mathbb R}^{2})$ 
without loosing the important compactness of weakly
singular integral operators.

At the core of the Mikhlin-Pr{\"o}ssdorf approach, sucessfully applied in \cite{bud_sam_06, SamokhinBook}
to the three-dimensional volume integral equation, is the symbol calculus. It allows to
study the existence of a solution and shows a way to reduce
the original singular integral operator to a manifestly Fredholm 
operator of the form `identity plus compact'. The symbol of an integral operator in the present 
two-dimensional case is a matrix-valued function $\Phi_{nm}(\bx,\bk)$, $n,m=1,2,3$,
$\bx\in{\mathbb R}^{2}$, and $\bk\in{\mathbb R}^{2}$. Construction of the symbol 
follows a few simple rules explained in \cite{Mikhlin, SamokhinBook, bud_sam_06}:
the symbol of the sum/product of operators is the sum/product of symbols;
the symbol of a compact operator is zero; the symbol of a multiplication
operator is the multiplier function itself; the symbol of the strongly 
singular integral operator is the principal value of the Fourier transform 
of its kernel (that is where the second argument $\bk$ -- the dual 
of $\bx$ -- comes from). It follows, in particular, that the symbol of the 
matrix operator is a matrix-valued function. Applying these rules 
to (\ref{eq:SingularTE}) we obtain
\begin{align}
\label{eq:Symbol}
{\mathbf \Phi}(\bx,\bk)=
\begin{bmatrix}
s(\bx) & 0 & 0 \\
0 & s(\bx) & 0 \\
0 & 0 & 1
\end{bmatrix}
+
\begin{bmatrix}
\tilde{G}_{11}(\bk) & \tilde{G}_{12}(\bk) & 0 \\
\tilde{G}_{21}(\bk) & \tilde{G}_{22}(\bk) & 0 \\
0 & 0 & 0
\end{bmatrix}
\begin{bmatrix}
\chi_{\rm e}(\bx) & 0 & 0 \\
0 & \chi_{\rm e}(\bx) & 0 \\
0 & 0 & \chi_{\rm m}(\bx)
\end{bmatrix}.
\end{align}
The Fourier transforms $\tilde{G}_{nm}(\bk)$ of 
the kernels $G_{nm}(\bx)$ are found as a series of two-dimensional `spherical' functions of order $p$:
\begin{align}
\label{eq:SphericalExpansion}
\tilde{G}_{nm}(\bk)=\sum\limits_{p=0}^{\infty}
\gamma_{2,p}\left[a^{(1)}_{p}Y^{(1)}_{p,2}\left(\frac{\bk}{\vert\bk\vert}\right)
+a^{(2)}_{p}Y^{(2)}_{p,2}\left(\frac{\bk}{\vert\bk\vert}\right)\right],
\end{align}
which, in general, results in a different expression for each combination of $n$ and $m$. 
The particular form of (\ref{eq:SphericalExpansion}) reflects the fact
that $\tilde{G}_{nm}$ is a function of the direction of $\bk$ only and it does 
not depend on its magnitude (see \cite{bud_sam_06, Mikhlin} for a formal proof of this fact).
In two dimensions, with $\tilde{\phi}$ denoting the
directional angle of the unit vector $\bk/\vert\bk\vert$ we have: 
$Y^{(1)}_{p,2}=\sin(p\tilde{\phi})$, $Y^{(2)}_{p,2}=\cos(p\tilde{\phi})$,
and $\gamma_{2,p}=\pi i^{p} \Gamma(p/2)/\Gamma((p+2)/2)$.
The expansion coefficients in
(\ref{eq:SphericalExpansion}) can be easily recovered by representing
the original kernels $G_{nm}(\bx)$ in the following form:
\begin{align}
\label{eq:BasicKernel}
G_{nm}\left(\bx\right)=
\frac{1}{\vert\bx\vert^{2}}f_{nm}\left(\frac{\bx}{\vert\bx\vert}\right).
\end{align}
If we now express the characteristic $f_{nm}$ as
\begin{align}
\label{eq:Characteristics}
f_{nm}\left(\frac{\bx}{\vert\bx\vert}\right)=
\sum\limits_{p=1}^{\infty}\left[a^{(1)}_{p}\sin(p\phi)+a^{(2)}_{p}\cos(p\phi)\right],
\end{align}
where $\phi$ is the directional angle of the unit vector $\bx/\vert\bx\vert$,
then the expansion coefficients in (\ref{eq:SphericalExpansion}) are those of
(\ref{eq:Characteristics}). Again we note that for our matrix-valued kernel 
this has to be done for each combination of $n$ and $m$ separately.
Eventually we arrive at
\begin{align}
\label{eq:SymbolSO}
\begin{bmatrix}
\tilde{G}_{11}(\bk) & \tilde{G}_{12}(\bk) & 0 \\
\tilde{G}_{21}(\bk) & \tilde{G}_{22}(\bk) & 0 \\
0 & 0 & 0
\end{bmatrix}
=
\begin{bmatrix}
\cos^{2}(\tilde{\phi}) - 1/2 & \sin(\tilde{\phi})\cos(\tilde{\phi}) & 0 \\
\sin(\tilde{\phi})\cos(\tilde{\phi}) & \sin^{2}(\tilde{\phi}) - 1/2 & 0 \\
0 & 0 & 0
\end{bmatrix}.
\end{align}
Finally the symbol of the complete operator is
\begin{align}
\label{eq:CompleteSymbol}
{\mathbf \Phi}(\bx,\tilde{\phi})= 
\begin{bmatrix}
1 & 0 & 0 \\
0 & 1 & 0 \\
0 & 0 & 1
\end{bmatrix}
+
\chi_{\rm e}(\bx)
\begin{bmatrix}
\cos^{2}(\tilde{\phi}) & \sin(\tilde{\phi})\cos(\tilde{\phi}) & 0 \\
\sin(\tilde{\phi})\cos(\tilde{\phi}) & \sin^{2}(\tilde{\phi}) & 0 \\
0 & 0 & 0
\end{bmatrix}
={\mathbb I}+\chi_{\rm e}(\bx){\mathbb Q},
\end{align}
where matrix ${\mathbb Q}$ is a projector, i.e., ${\mathbb Q}^{2}={\mathbb Q}$.

Several important conclusions can be made already at this stage.
First of all we notice an almost complete equivalence of the symbol
in the two-dimensional TE case with the previously obtained 
symbol of the three-dimensional scattering operator 
\cite{bud_sam_06, SamokhinBook}. The necessary and sufficient
condition for the existence of a solution is now readily obtained 
as the condition on the invertibility of ${\mathbf \Phi}$, which is
\begin{align}
\label{eq:ExistenceTE}
\varepsilon(\bx,\omega)\neq 0, \;\;\;\bx\in{\mathbb R}^{2}.
\end{align}
In the TM case (\ref{eq:MatrixSolutionTM}) the symbol can be 
derived along the same lines and turns out to 
be ${\mathbb I}+\chi_{\rm m}(\bx){\mathbb Q}$, with the existence 
condition being $\mu(\bx,\omega)\neq 0$, $\bx\in{\mathbb R}^{2}$.
It is interesting to note the difference in these existence conditions 
with respect to the three-dimensional case, where both $\varepsilon(\bx)$
and $\mu(\bx)$ are required to be nonzero at the same time.

The explicit inverse of the symbol matrix is found to be
\begin{align}
\label{eq:InverseSymbol}
{\mathbf \Phi}^{-1}(\bx,\tilde{\phi})={\mathbb I}+\chi'_{\rm e}(\bx){\mathbb Q},
\end{align}
where $\chi'_{\rm e}(\bx,\omega)=\varepsilon_{\rm b}/\varepsilon(\bx,\omega)-1$.
The inverse of the symbol is the symbol of the regularizer -- 
an operator, which reduces the original singular integral operator
to the form `identity plus compact'. 
Since ${\mathbf \Phi}^{-1}$
has the same form as ${\mathbf \Phi}$, we conclude that 
our original
operator with its electric contrast function changed into $\chi'_{\rm e}$ 
may be employed as a regularizer.
Thus, if $A(\chi_{\rm e})$ is the original operator from 
(\ref{eq:StandardForm})--(\ref{eq:SingularTE}),
and $A(\chi'_{\rm e})$ denotes the same operator with the contrast 
function replaced by $\chi'_{\rm e}$, then we have
\begin{align}
\label{eq:Regularization}
A(\chi'_{\rm e})A(\chi_{\rm e})=I+K,
\end{align}
where $K$ is a generic compact operator.
This looks strikingly similar to the Calder{\'o}n identity \cite{Andriulli_2008, Bagci_2010, Christiansen_2003},
and can be viewed as an extension of the latter to the case of an inhomogeneous penetrable scatterer.

Regularizer eliminates the essential spectrum of $A(\chi_{\rm e})$. 
Indeed, even if $A(\chi_{\rm e})$ has a nontrivial essential 
spectrum, densely spread over a part of the complex plane,
the spectrum of $I+K$ consists of isolated 
points only. As can be deduced from the symbol (\ref{eq:CompleteSymbol}),
the essential spectrum of the operator $A(\chi_{\rm e})$ in the TE 
case is given by
\begin{align}
\label{eq:EssentialSpectrum}
\lambda_{\rm ess}=\frac{\varepsilon(\bx)}{\varepsilon_{\rm b}}, 
\;\;\;\bx\in{\mathbb R}^{2},
\end{align}
and by $\lambda_{\rm ess}=\mu(\bx)/\mu_{\rm b}$, 
$\bx\in{\mathbb R}^{2}$ in the TM case.
From the condition of the uniqueness of the solution we also
obtain a traditional \cite{Kleinman_Roach_VanDenBerg_1990} wedge-shaped bound on the possible discrete eigenvalues:
\begin{align}
\label{eq:eig_bound}
\begin{split}
&{\rm Im\,}\varepsilon(\bx)-\left[{\rm Im\,}\varepsilon(\bx)+{\rm Im\,}\varepsilon_{\rm b}\right]{\rm Re\,}\lambda\\
&+\left[{\rm Re\,}\varepsilon(\bx)-{\rm Re\,}\varepsilon_{\rm b}\right]{\rm Im\,}\lambda
+{\rm Im\,}\varepsilon_{\rm b}\vert\lambda\vert^2\leq0,\quad\bx\in D.
\end{split}
\end{align}

\section{GMRES convergence and matrix eigenvalues}
Let us introduce a uniform two-dimensional $N$-node grid over a 
rectangular computational domain $D$. For simplicity let the grid step 
$h$ be the same in both directions. 
Typically the grid step is chosen in accordance with the following 
empirical rule: find out the highest local value of the refractive 
index $n={\rm max}\{\sqrt{\varepsilon'(\bx)\mu'(\bx)]/
\varepsilon_{\rm b}\mu_{\rm b}}\}$, $\bx\in D$, where prime denotes 
the real part; choose an integer $k$ -- the number of points 
per smallest wavelength; and, finally, compute the grid step as 
$h=\lambda_{\rm b}/(kn)$, where $\lambda_{\rm b}$ is the wavelength in the background medium. 
We set our discretization here
at the traditional level of $k=15$ points, which,
as we verified, gives a small global error (in the order of $10^{-5}$) with respect 
to the analytical solution for a homogeneous circular cylinder. For rectangular grid-conforming 
boundaries the error with this discretization rule could be even smaller, since most of the error
for a circular cylinder comes from a poor geometrical representation of the boundary
and/or bad approximation of the area of the circular cross-section.

Applying a simple collocation
technique with the mid-point rule we end up with an algebraic 
system of the order $3N\times 3N$:
\begin{align}
 \label{eq:OurSystem}
Au=b,
\end{align}
with the following structure:
\begin{align}\label{eq:lin_system}
\begin{bmatrix}
A_{11} & A_{12} & A_{13} \\
A_{21} & A_{22} & A_{23} \\
A_{31} & A_{32} & A_{33}
\end{bmatrix}
\begin{bmatrix}
u_1 \\
u_2 \\
u_3
\end{bmatrix}=
\begin{bmatrix}
b_1 \\
b_2 \\
b_3
\end{bmatrix},
\end{align}
where the matrix elements are
\begin{align}\label{eq:A11}
\begin{split}
[A_{\ell\ell}]_{nm}=&-k_{\rm b}^2h^2\chi_{\rm e}(\bx_m)\bigg\{\left[\frac{i}{2k_{\rm b}r_{nm}}H_1^{(1)}(k_{\rm b}r_{nm})-\frac{i}{4}H_0^{(1)}(k_{\rm b}r_{nm})\right]\theta_{\ell,nm}\theta_{\ell,nm}\\
&+\left[\frac{i}{4}H_0^{(1)}(k_{\rm b}r_{nm})-\frac{i}{4k_{\rm b}r_{nm}}H_1^{(1)}(k_{\rm b}r_{nm})\right]\delta_{\ell\ell}\bigg\},\quad \ell=1,2,\quad m\neq n
\end{split}
\end{align}
\begin{align}
\begin{split}
[A_{\ell\ell}]_{nn}=1+\left[1-\frac{i\pi k_{\rm b}h}{4\sqrt{\pi}}H_1^{(1)}(k_{\rm b}h/\sqrt{\pi})\right]\chi_{\rm e}(\bx_n),\quad \ell=1,2
\end{split}
\end{align}
\begin{align}
\begin{split}
[A_{\ell q}]_{nm}=&-k_{\rm b}^2h^2\chi_{\rm e}(\bx_m)\left[\frac{i}{2k_{\rm b}r_{nm}}H_1^{(1)}(k_{\rm b}r_{nm})-\frac{i}{4}H_0^{(1)}(k_{\rm b}r_{nm})\right]\\
&\times\theta_{\ell,nm}\theta_{q,nm}(1-\delta_{nm}),\quad\ell,q=1,2,\quad\ell\neq q
\end{split}
\end{align}
\begin{align}
\begin{split}
[A_{13}]_{nm}=i\omega\mu_{\rm b}h^2\chi_{\rm m}(\bx_m)\frac{ik_{\rm b}}{4}H_1^{(1)}(k_{\rm b}r_{nm})\theta_{2,nm}(1-\delta_{nm})
\end{split}
\end{align}
\begin{align}
\begin{split}
[A_{32}]_{nm}=-i\omega\varepsilon_{\rm b}h^2\chi_{\rm e}(\bx_m)\frac{ik_{\rm b}}{4}H_1^{(1)}(k_{\rm b}r_{nm})\theta_{1,nm}(1-\delta_{nm})
\end{split}
\end{align}
\begin{align}
\begin{split}
[A_{31}]_{nm}=-[A_{13}]_{nm}
\end{split}
\end{align}
\begin{align}
\begin{split}
[A_{23}]_{nm}=-[A_{32}]_{nm}
\end{split}
\end{align}
\begin{align}
\begin{split}
[A_{33}]_{nm}=-k_{\rm b}^2h^2\chi_{\rm m}(\bx_m)\frac{i}{4}H_0^{(1)}(k_{\rm b}r_{nm}),\quad m\neq n
\end{split}
\end{align}
\begin{align}\label{eq:A33}
\begin{split}
[A_{33}]_{nn}=1+\left[1-\frac{i\pi k_{\rm b}h}{2\sqrt{\pi}}H_1^{(1)}(k_{\rm b}h/\sqrt{\pi})\right]\chi_{\rm m}(\bx_n)
\end{split}
\end{align}
In these formulas: $r_{nm}=\vert\bx_{n}-\bx_{m}\vert$, $
\theta_{\ell,nm}=(x_{\ell,n}-x_{\ell,m})/r_{nm}$, $\ell=1,2$, $n,m=1,\ldots,N$. 
Here $x_{\ell,n}$ denotes the Cartesian component of the two-dimensional
position vector $\bx_{n}$ pointing at the $n$th node of the grid.
The unkhown vector $u=[u_{1},u_{2},u_{3}]^T$ contains the grid values of the TE total 
field $[E_1(\bx_{m}),E_2(\bx_{m}),H_3(\bx_{m})]^T$, 
$m=1,\dots,N$. The right-hand side vector $b=[b_{1},b_{2},b_{3}]^T$
contains the grid values of the incident field. For example, the TE plane wave impinging to the object 
at an angle $\psi$ with respect to $x_1$-axis is represented by: 
$[b_1]_{n}=\exp[ik_{\rm b}(x_{1,n}\cos\psi+x_{2,n}\sin\psi)]$, 
$[b_2]_{n}=-[b_1]_{n}k_{\rm b}\sin\psi/(\omega\varepsilon_{\rm b})$, and 
$[b_3]_{n}=[b_1]_{n}k_{\rm b}\cos\psi/(\omega\varepsilon_{\rm b})$, $n=1,\dots,N$.

\begin{figure}[t!]
\centering
   \includegraphics[scale=1]{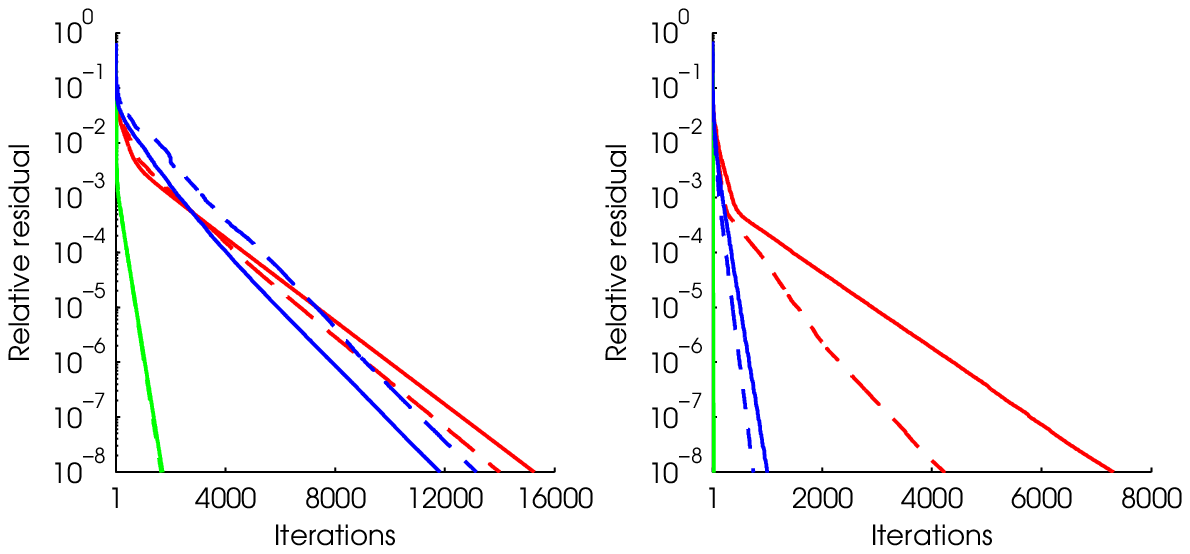}
\caption{Convergence of restarted GMRES applied to the original system 
with (dashed) and without (solid) deflation of the largest 
eigenvalues for three physically distinct scatterers. Deflation accelerates GMRES 
in the TM case (right), but does not work in the TE case (left). 
lossless object -- red lines, object with small losses and negative real part of permittivity -- green lines, 
and inhomogeneous object -- blue lines.}
\label{precond_original}
\end{figure}
Both the system matrix and the original integral operator 
are neither symmetric nor normal. Therefore, the range of iterative methods applicable to the problem
is extremely limited, with (restarted) GMRES typically showing the best performance.
\fig{precond_original} shows the GMRES convergence histories on the original system $Au=b$ for three
physically distinct scatterers (all with $\chi_{\rm m}=0$) and two different polarizations, corresponding to the 
TE and TM cases. The three scatterers are: a lossless object with large positive real
permittivity, $\varepsilon/\varepsilon_{\rm b}=16$, an object with small losses and a large negative real part of the 
permittivity, $\varepsilon/\varepsilon_{\rm b}=-16+i1.5$, 
and an inhomogeneous object consisting of two concentric layers with the inner part having large losses, 
$\varepsilon/\varepsilon_{\rm b}=2.5+i20$, and the outer layer having a large positive real permittivity, 
$\varepsilon/\varepsilon_{\rm b}=16$. All objects are cylinders with square cross-sections
with the side length $a=\lambda_{\rm b}$. The inner core of the inhomogeneous object has also a square 
cross-section with the side length 
$a/2$. The {\tt restart} parameter of the GMRES algorithm is equal to $40$ (we shall study the influence of this
parameter in more detail later). 
In \fig{precond_original} (left) we illuminate the scatterers with a plane wave whose electric field vector is 
orthogonal to the axis of the cylinder (TE), whereas in \fig{precond_original}(right)
the electric field vector of the incident plane wave is parallel to the axis of the cylinder (TM).
Since all three objects do not exhibit any magnetic contrast with respect to the background medium, the
integral operators are different for these TE and TM cases -- see (\ref{eq:MatrixSolutionTE}), (\ref{eq:MatrixSolutionTM}). 
The latter contains now only the weakly singular part and is, therefore, of the form `identity plus compact'. 

Next to the convergence plots of the original system in \fig{precond_original} 
we present the convergence histories of the restarted  
GMRES with a deflation-based right preconditioner (dashed lines). As was suggested in \cite{sifuentes}, 
deflating the largest eigenvalues of $A$ in the TM case may significantly accelerate the convergence of the full (un-restarted) GMRES. 
The preconditioner deflating $r$ largest eigenvalues of $A$ is constructed by first running the {\tt eigs}
algorithm to retrieve these $r$ eigenvalues and the corresponding eigenvectors. With the help of the modified
Gramm-Schmidt algorithm we further build an orthonormal basis for the retrieved set of eigenvectors.
Let this basis be stored in a $3N\times r$ matrix $V$. Then, the right preconditioner \cite{erh_bur_poh_96} has the form
\begin{align}
 \label{eq:DeflatingPrecondtioner}
 P^{-1}=I_{3N}+V\left[T^{-1}-I_r\right]V^\ast,
\end{align}
with
\begin{align}
T=V^\ast AV.
\end{align}
As was observed in \cite{sifuentes} for the TM case, although the system matrix is not normal, 
the matrix $V$ has full column rank, and the matrix $T$ has a decent 
condition number. We confirm this behavior in the present TE case. 
The term `deflation' refers to the fact that the spectrum of $AP^{-1}$ looks like the 
spectrum of $A$ with $r$ largest eigenvalues moved to the point $1+i0$ of the complex plane,
and the corresponding eigen-subspace has been projected out. 
Here, even with the restarted GMRES, we see some improvement in the TM case \fig{precond_original}(right).
The TE case, however, is not affected by deflation (it may even become worse).  

\tab{tab:conv_his_TM} and \tab{tab:conv_his_TE} summarize 
the numerical results corresponding to the test objects in the TM  and TE cases.
\begin{table}[t!]
\caption{Restarted GMRES with fair memory usage for three different test scatterers, TM case. Original and deflated systems.}
\begin{center} \footnotesize
\begin{tabular}{|c|c|c|c|c|c|c|} \hline  
               &                                                              & $r$, deflated   &                         & CPU time,   & Iterations,          & Speed- \\ 
System & $\varepsilon/\varepsilon_{\rm b}$ & eigenvalues    & {\tt restart} & seconds      & tol. $10^{-8}$  & up         \\ \hline\hline
$A$ & 16 & 0 & 40 & 168 & 7311 & 1.0 \\ \hline
$AP^{-1}$ & 16 & 28 & 12 & 78 & 4230 & 2.2 \\ \hline
$A$ & $-16+i1.5$ & 0 & 40 & 0.7 & 27 & 1.0 \\ \hline
$AP^{-1}$ & $-16+i1.5$ & 28 & 12 & 0.5 & 14 & 1.4 \\ \hline
$A$ & $16$ and $2.5+i20$ & 0 & 40 & 23 & 1000 & 1.0 \\ \hline
$AP^{-1}$ & $16$ and $2.5+i20$ & 28 & 12 & 14 & 749 & 1.6 \\ \hline
\end{tabular}
\end{center} 
\label{tab:conv_his_TM}
\end{table}
The $3N\times r$ matrix $V$ needs to be stored in the memory. For the sake of fair comparison,
we divide our memory between the inner Krylov subspace of the restarted GMRES and $V$, 
i.e., if we set ${\tt restart}=k$, $k>r$ 
with the original matrix $A$, then we use ${\tt restart}=k-r$ with the preconditioned matrix $AP^{-1}$.

\begin{figure}[t!]
\centering
\subfigure[TE case.]{
   \includegraphics[scale=1]{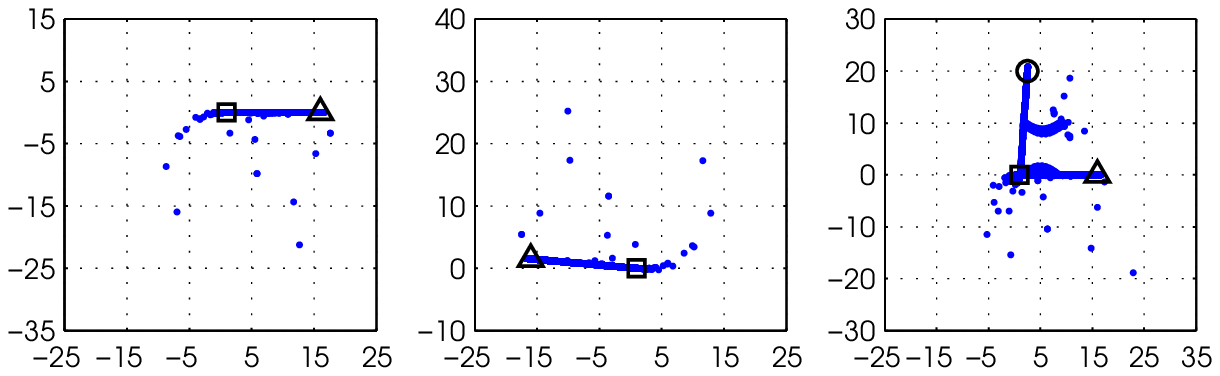}
   \label{spectrum_TE}
 }
 \subfigure[TM case.]{
   \includegraphics[scale=1]{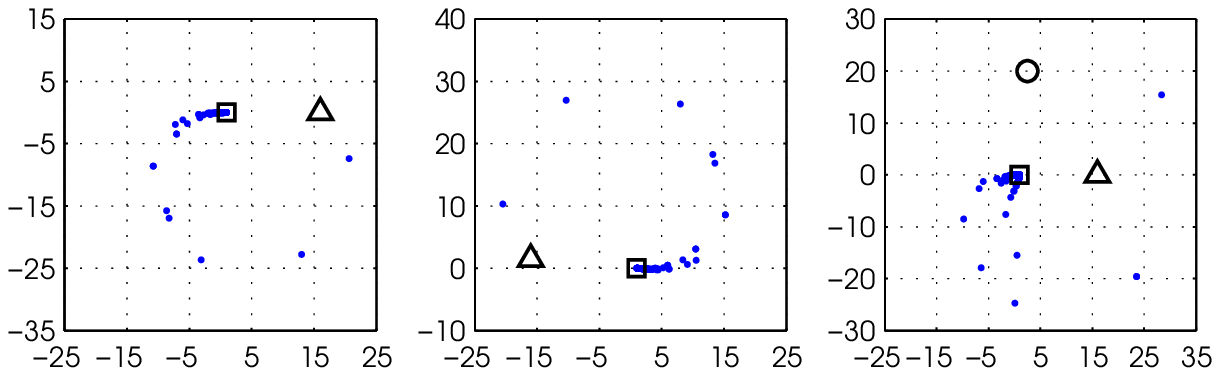}
   \label{spectrum_TM}
 }
\caption{Eigenvalues of the system matrix $A$ for three physically distinct 
scatterers in the TE (top) and TM (bottom) cases. 
Horizontal axis -- real part, vertical axis -- imaginary part.
Also shown (squares, triangles and circles) the values of the relative permittivity, including the background medium,
in each scattering configuration.}
\label{spectra_original}
\end{figure}

One can understand the peculiar lack of progress in the TE case by considering  \fig{spectra_original},
which shows the spectra of the 
system matrix for the TE and TM cases corresponding to the numerical experiments of \fig{precond_original}. 
As can be seen, the difference between the spectra is in
the presence of dence line segments $[1,\varepsilon(\bx)/\varepsilon_{\rm b}]$, $\bx\in{\mathbb R}^{2}$  
and some extra isolated eigenvalues inside a distribution remotely resembling a shifted circumference 
(the latter is much more pronounced with homogeneous objects at higher frequencies, see \cite{sifuentes}). 
The dense line segments in the TE spectrum contain a very large number of closely spaced eigenvalues,
and the deflation process simply stagnates when it reaches the outer ends of these lines,
corresponding to $\varepsilon(\bx)/\varepsilon_{\rm b}$, $\bx\in D$.
Indeed, in the TE case we have seen no change in the convergence with $AP^{-1}$ when we tried to deflate
any eigenvalues beyond those that were present outside a circle with the center at zero and the radius
reaching the largest value of $\vert\varepsilon(\bx)/\varepsilon_{\rm b}\vert$, $\bx\in{\mathbb R}^{2}$
for each particular scattering configuration. At the same time 
in the TM case we observed consistent improvement in convergence when deflating 
all eigenvalues outside the unit circle, see \tab{tab:conv_his_TM}. 

\section{Preconditioning}
In the three-dimensional case it was conjectured \cite{bud_sam_06} that the dense line segments, 
which we now observe in the TE spectrum as well, are the discrete image of the essential spectrum,
$\lambda_{\rm ess}=\varepsilon(\bx)/\varepsilon_{\rm b}$, of the corresponding singular 
integral operator. The absence of these lines in the TM spectrum is just another indirect confirmation.
Yet, we do not have a formal proof, since there is no such thing as essential spectrum with matrices,
and its role and transformation in the process of discretization is unclear. Anyway, the process of
deflation is obviously hindered by those dense line segments and, assuming that they are due to the
essential spectrum, we already know how to `deflate' them. 
On the continuous level this can be achieved by regularization.
We shall, therefore, apply the discretized regularizer $A_{\rm R}$, i.e., a discrete version of the  
operator $A(\chi'_{\rm e})$,
so that the problem becomes $A_{\rm R}Au=A_{\rm R}b$. Notice that no additional storage is required,
since the original Green's tensor can be re-used. Matrix-vector products, however, do become twice 
more expensive to compute.

\begin{figure}[t!]
\centering
\includegraphics[scale=1]{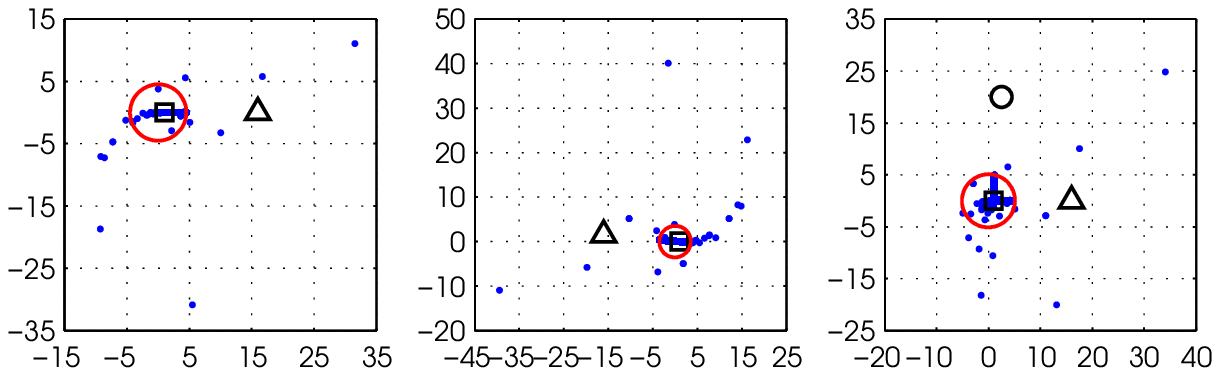}
\caption{Eigenvalues of the regularized matrix $A_{\rm R}A$, corresponding to the spectra of $A$
in the TE case, see \fig{spectra_original} (top). Red circles show the estimated outer bounds of 
the dense cluster -- deflation bounds.}
\label{spectra_reg}
\end{figure}

To see how this procedure affects the matrix spectrum we compute the eigenvalues of $A_{\rm R}A$
corresponding to the matrices analyzed in \fig{spectra_original} (top). The results are presented in 
\fig{spectra_reg}, where we see that the line segments have been substantially reduced.
For a homogeneous object we can see what is going on also on the matrix level. The general form of our system 
matrix is $A=I-GX$, where $I$ is a $3N\times 3N$ identity matrix, $X$ is the diagonal matrix, containing
the grid values of the contrast function, and $G$ is a dense matrix produced by the integral operator 
with the Green's tensor kernel. When the contrast is homogeneous and equal to a scalar complex value 
$\chi=\varepsilon/\varepsilon_{\rm b}-1$ inside the scatterer, we can write our matrix as $A=I-\chi G$. Now, let $\psi_{n}$ be an eigenvector
of $G$, i.e. $G\psi_{n}=\beta_{n}\psi_{n}$. Then, the eigenvalues of $A$ are:
\begin{align}
 \label{eq:EigenvaluesA}
 \lambda_{n}=1-\beta_{n}\left(\frac{\varepsilon}{\varepsilon_{\rm b}}-1\right),
\end{align}
and the eigenvalues of $A_{\rm R}A$ are
\begin{align}
 \label{eq:EigenvaluesArAbeta}
 \lambda'_{n}=\left[1-\beta_{n}\left(\frac{\varepsilon_{\rm b}}{\varepsilon}-1\right)\right]\lambda_{n}.
\end{align}
Eliminating $\beta_{n}$ we get
\begin{align}
 \label{eq:EigenvaluesArA}
 \lambda'_{n}=\left[1+\frac{\varepsilon_{\rm b}}{\varepsilon}(1-\lambda_{n})\right]\lambda_{n}.
\end{align}
We have verified numerically that this transformation holds for the first two plots of
\fig{spectra_original} (top) and \fig{spectra_reg}. Curiously, the seemingly quadratic
amplification of the large eigenvalues in (\ref{eq:EigenvaluesArA}) is moderated by 
the inverse of the relative permittivity. That is why the regularized spectra are not as 
extended as one would fear.
Unfortunately, such a simple proof is not possible with an
inhomogeneous object, yet the effect of the regularization on the matrix spectrum in 
\fig{spectra_original} (top, right) and \fig{spectra_reg} (right) is obviously similar to the 
homogeneous cases.

\tab{tab:conv_his_TE} shows that such a regularized system is already better than the original one in 
terms of convergence. 
The relative speed-up, however, is limited to approximately $2.5$ times 
due to the extra work involved in computing the matrix-vector products with $A_{\rm R}A$.
Moreover, no convergence improvement is observed  with the negative permittivity object.

To further accelerate the convergence and make the method work with the negative permittivity 
objects as well (they are important in surface plasmon-polariton studies) we can now apply the deflation technique. 
For that we need to decide on the number of eigenvalues to be deflated. Obviously, it does not
make sense to deflate beyond the first (from outside) dense cluster -- deflation will only 
cost more time and memory, while the improvement will stagnate.
After the regularization, the outer boundaries of the dense line segments have been all shifted 
inwards. However, due to nonlinear nature of the mapping (\ref{eq:EigenvaluesArA})
it is hard to tell exactly what an image of an arbitrary line segment would look like. 
Only for a lossless homogeneous scatterer are we able to predict that a real line segment
$[1,\varepsilon/\varepsilon_{\rm b}]$ will be mapped onto the real line segment
$[1,1+(1-\varepsilon/\varepsilon_{\rm b})^{2}/(4\varepsilon/\varepsilon_{\rm b})]$,
which contracts the original segment by approximately four times.
For a general line segment radially emerging from the point $1+i0$ we have
made a small numerical routine, which produces a discretized image according to (\ref{eq:EigenvaluesArA})
and finds the point with the largest absolute value -- the radius of the circle beyond which
all eigenvalues should be deflated. We have observed, by analyzing the spectra for
various inhomogeneous objects, that the general regularization seems to transform
each line segment according to (\ref{eq:EigenvaluesArA}). Therefore, in an 
inhomogeneous case, we can simply apply the aforementioned
numerical mapping to each of the segments $[1,\varepsilon(\bx)/\varepsilon_{b}]$, $\bx\in D$,
and choose the point with the maximum absolute value as the deflation radius.
The deflation bounds obtained in this way are shown in \fig{spectra_reg} as solid red circles.

\begin{figure}[t!]
\centering
   \includegraphics[scale=1]{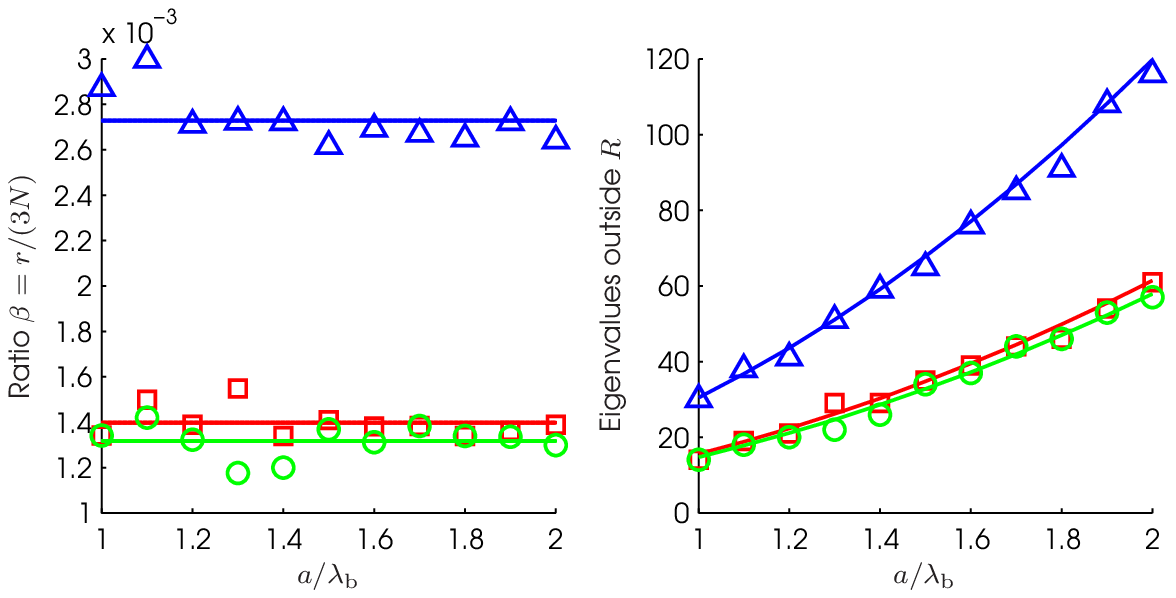}
\caption{Dependence of the number of eigenvalues to be deflated on the wavelength
for three test scatterers (same medium parameters as in previous figures). 
Left: ratio $\beta$ of the number $r$ of eigenvalues outside the deflation bound 
(circles with radius $R$, see \fig{spectra_reg}) to the total number $3N(\lambda_{\rm b})$ of eigenvalues for 
a classical discretization rule (fixed number of points per smallest medium wavelength). 
Right: actual number $r$ of eigenvalues to be deflated (with magnitude greater than the 
deflation bound $R$) and our analytical prediction based on an average $\beta$ coefficient 
estimated from the left plot.}
\label{ratios}
\end{figure}

The number of eigenvalues outside the deflation bound depends on the physics of the 
problem. For example, we know that at low frequencies (large wavelengths) all discrete 
eigenvalues are located within the convex hull of the essential spectrum \cite{Budko_2005},
and that at higher frequencies they tend to spread out. We have investigated this phenomenon
in more detail and have found a pattern, similar to the one reported in \cite{sifuentes}.
As one can see from \fig{ratios}~(left), the ratio of the number $r$ of eigenvalues outside the
deflation bound to the total number of eigenvalues $3N$ is roughly constant (for objects larger than $\lambda_{\rm b}$),
if at each wavelength we adjust the spatial discretization according to the classical rule (say, $15$ points per
smallest medium wavelength). Hence, with fixed medium parameters and a square 
computational domain with side $a$, the total number 
of unknowns and therefore eigenvalues grows as 
$3N=3[\alpha^2(a/\lambda_{\rm b})^2+2\alpha(a/\lambda_{\rm b})+1]$ with 
$\alpha=15\max\{\sqrt{(\varepsilon'/\varepsilon_{\rm b})}\}$. Thus, knowing $\beta=r/(3N)$ 
for some wavelength, we can estimate $r=3N\beta$ at any other wavelength. 
This empirical law comes very handy since solving a relatively low-frequency spectral problem can 
be really easy due to small system matrix dimensions. Otherwise, if there is no shortage of memory for
the chosen wavelength, one can simply run the {\tt eigs} routine increasing the number of 
recovered largest eigenvalues until the deflation bound obtained above is reached.

Having a good estimate of the deflation parameter $r$ we have applied the right preconditioner
given by (\ref{eq:DeflatingPrecondtioner}) to the regularized system for the same three objects
as in all previous examples. The largest eigenvalues/eigenvectors of the regularized
system turn out to be very stable, so that one can safely apply the {\tt eigs} routine
with its tolerance set as high as $10^{-2}$. Though, we have not seen any substantial 
acceleration beyond the $10^{-4}$ tolerance, which, in its turn, gave us a 
significant speed-up in the off-line work compared to the default machine precision tolerance.
In all our experiments, the off-line time remained a small fraction of the total CPU time (for concrete 
numbers we refer to the example at the end of this Section).

\begin{figure}[t!]
\centering
   \includegraphics[scale=1]{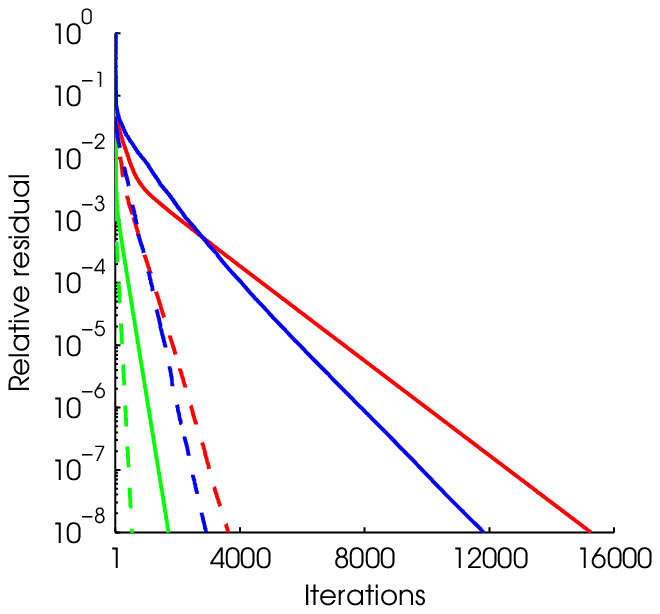}
\caption{Convergence of restarted GMRES with the original $Au=b$ (solid) and preconditioned 
$A_{\rm R}AP^{-1}v=A_{\rm R}b$ (dashed) systems for the three physically distinct scatterers in the TE case. 
Red lines; $\varepsilon/\varepsilon_{\rm b}=16$; Green lines: $\varepsilon/\varepsilon_{\rm b}=-16+i1.5$;
Blue lines: inhomogeneous object $\varepsilon/\varepsilon_{\rm b}=(16,2.5+i20)$.}
\label{conv_RD}
\end{figure}

The convergence histories for a restarted GMRES with this preconditioner 
are presented in \fig{conv_RD}. 
\begin{table}[t!]
\caption{Restarted GMRES with fair memory usage for three different test scatterers, TE case. Original, deflated, 
regularized, and regularized-plus-delfated systems.}
\begin{center} \footnotesize
\begin{tabular}{|c|c|c|c|c|c|c|} \hline  
               &                                                              & $r$, deflated   &                         & CPU time,   & Iterations,          & Speed- \\ 
System & $\varepsilon/\varepsilon_{\rm b}$ & eigenvalues    & {\tt restart} & seconds      & tol. $10^{-8}$  & up         \\ \hline\hline
$A$ & 16 & 0 & 40 & 350 & 15258 & 1.0 \\ \hline
$AP^{-1}$ & 16 & 6 & 34 & 306 & 14002 & 1.1 \\ \hline
$A_{\rm R}A$ & 16 & 0 & 40 & 142 & 4030 & 2.5 \\ \hline
$A_{\rm R}AP^{-1}$ & 16 & 14 & 26 & 117 & 3626 & 3 \\ \hline
$A$ & $-16+i1.5$ & 0 & 40 & 39 & 1708 & 1.0 \\ \hline
$AP^{-1}$ & $-16+i1.5$ & 7 & 33 & 36 & 1659 & 1.1 \\ \hline
$A_{\rm R}A$ & $-16+i1.5$ & 0 & 40 & 36 & 1021 & 1.1 \\ \hline
$A_{\rm R}AP^{-1}$ & $-16+i1.5$ & 30 & 10 & 16 & 534 & 2.4 \\ \hline
$A$ & $16$ and $2.5+i20$ & 0 & 40 & 275 & 11823 & 1.0 \\ \hline
$AP^{-1}$ & $16$ and $2.5+i20$ & 3 & 37 & 310 & 13156 & 0.9 \\ \hline
$A_{\rm R}A$ & $16$ and $2.5+i20$ & 0 & 40 & 107 & 3032 & 2.6 \\ \hline
$A_{\rm R}AP^{-1}$ & $16$ and $2.5+i20$ & 14 & 26 & 94 & 2911 & 3 \\ \hline
\end{tabular}
\end{center} 
\label{tab:conv_his_TE}
\end{table}
This figure, and the iteration counts in general, are slightly misleading
when one deals with an iterative algorithm like restarted GMRES. For example, we see in \tab{tab:conv_his_TE}
that, apart from the negative permittivity case, the iteration counts after deflation 
are only slightly better than what we have already
achieved with the regularized system. Hence, one could conclude that the deflation 
is not worth the effort. However, the GMRES algorithm does not scale lineraly in time as a function 
of the {\tt restart} parameter, i.e., while it may converge in less iterations for larger values of {\tt restart},
the CPU time may not decrease as much. This has partly to do with the re-orthogonalization process, 
which takes more time for larger inner Krylov subspaces, and partly with the very bad spectral 
properties of our system matrix. 

Thus, to decide upon the benefits of deflation we compare execution times for various values of {\tt restart} 
parameter. In addition, we need to consider another important constraint -- the available memory. 
Since a fair comparison can only be achieved if we give the same amount of memory to
both the original and the preconditioned systems, we give 
all the available memory to the GMRES algorithm in the former case, 
and split it between the {\tt restart} and the deflation subspace of size $r$ in the latter.
Finally, we want to know if the efficiency of the proposed preconditioner depends on the physics of the 
problem, i.e., does it work better or worth for larger objects and contrasts?
Of course, the system becomes larger and more difficult to solve for larger objects, and we have to
be prepared to increase the value of {\tt restart} for larger objects and contrasts to achieve any 
convergence at all. Hence, assuming that we have enough memory at least
for deflation, we shall distinguish between the following three situations: 
limited memory, ${\tt restart}(A_{\rm R}AP^{-1})<r$;
moderate memory, ${\tt restart}(A_{\rm R}AP^{-1}) \sim r$; 
large memory, ${\tt restart}(A_{\rm R}AP^{-1}) > r$.

\begin{figure}[t!]
\centering
   \includegraphics[scale=1]{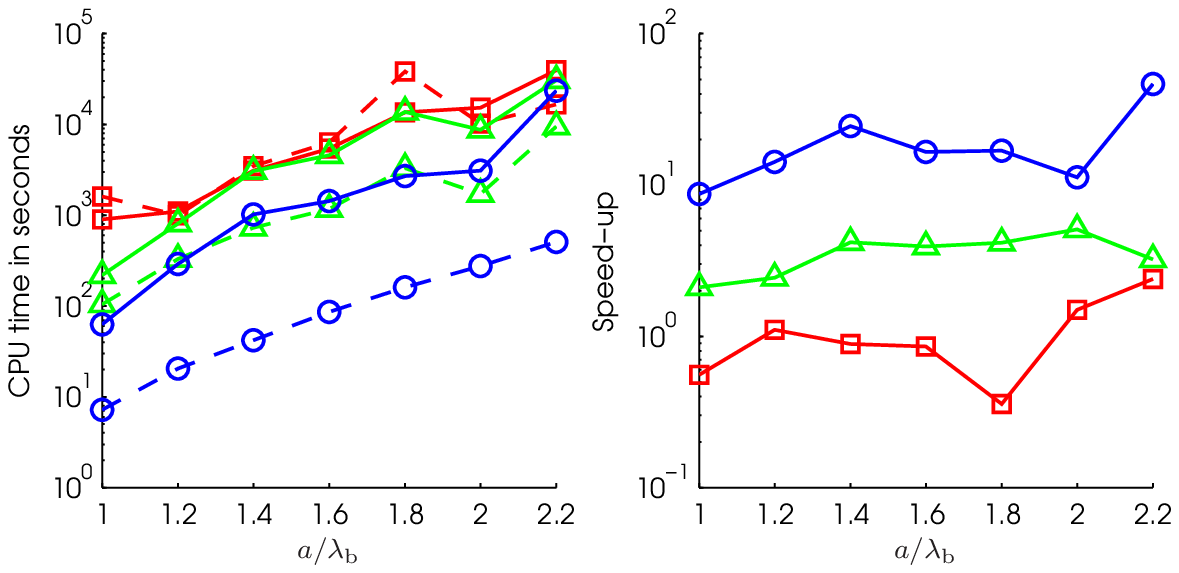}
\caption{CPU time (left) and the achieved speed-up (right) as a function of wavelength (object size)
for ${\tt restart}=r/2$ (red, squares), ${\tt restart}=2r$ (green, triangles), 
and ${\tt restart}=10r$ (blue, circles). 
The permittivty of the scatterer is $\varepsilon/\varepsilon_{\rm b}=16$.
Solid lines -- original system $Au=b$; dashed lines -- preconditioned system $A_{\rm R}AP^{-1}v=A_{\rm R}b$.}
\label{timing_16}
\end{figure}

We observe from the plots of  \fig{timing_16} that with limited memory, where we choose
${\tt restart}(A_{\rm R}AP^{-1})=r/2$ and ${\tt restart}(A)=r$,
the original system initially outperforms the preconditioned one. As the object size grows, so does
the size $r$ of the deflation subspace. Since in that case the ${\tt restart}(A_{\rm R}AP^{-1})$ parameter grows as well,
at some point (for $a/\lambda=1.9$ with this particular object) it reaches a value at which the preconditioned system
breaks even and begins to outperform the original one for larger objects (maximum achieved speed-up is $2.4$ times 
at $a/\lambda_{\rm b}=2.2$). 
With moderate memory, 
${\tt restart}(A_{\rm R}AP^{-1})=2r$, ${\tt restart}(A)=3r$, the preconditined system always performs 
better than the original one. Moreover, the relative speed-up \fig{timing_16}~(right) generally increases with the object size
(maximum achieved speed-up is $5.1$ times at $a/\lambda_{\rm b}=2$, minimum -- $2.1$ times at $a/\lambda_{\rm b}=1$).
With large memory, i.e., for ${\tt restart}(A_{\rm R}AP^{-1})=10r$ and ${\tt restart}(A)=11r$, the speed-up may be as 
high as $46.6$ times ($a/\lambda_{\rm b}=2.2$) (minimum speed-up of $8.7$ times at $a/\lambda_{\rm b}=1$), and it also 
grows (on average) with the object size. 

\begin{figure}[t!]
\centering
   \includegraphics[scale=1]{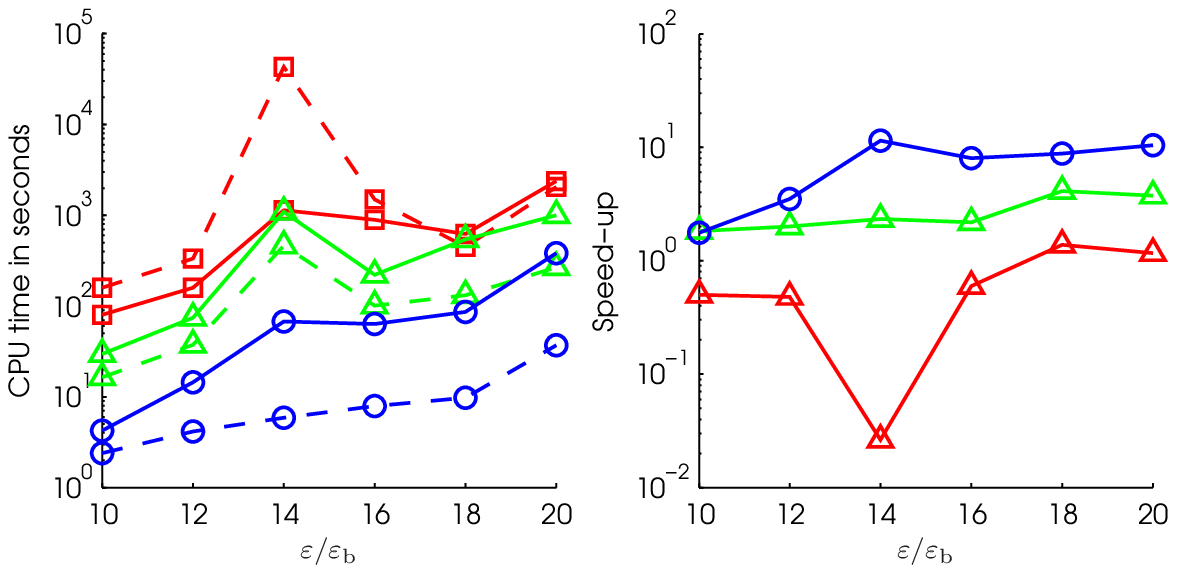}
\caption{CPU time (left) and the achieved speed-up (right) as a function of permittivity 
for ${\tt restart}=r/2$ (red, squares), ${\tt restart}=2r$ (green, triangles), 
and ${\tt restart}=10r$ (blue, circles). 
The scatterer is lossless with side length fixed at $a/\lambda_{\rm b}=1$.
Solid lines -- original system $Au=b$; dashed lines -- preconditioned system $A_{\rm R}AP^{-1}v=A_{\rm R}b$.}
\label{timing_16_er}
\end{figure}

With objects of finite extent changing the wavelength of the incident field is 
not the same as changing the object permittivity.
Hence, we repeat the same numerical experiments as in  \fig{timing_16_er}, but now for 
an object of fixed size, $a/\lambda_{\rm b}=1$, and varying permittivity. The deflation parameter $r$ is
estimated using the same algorithm as above and turns out to grow only slightly as a function of permittivity,
climbing from $13$ to $17$ for the permittivity changing between $10$ and $20$.
At the same time the total number of unknowns as well as eigenvalues grows very fast.
Therefore, deflation did not accelerate the iterative process here as much as when 
we changed the wavelength, and most of the speed-up must be due to the regularization.
For this relatively small scatterer we see some improvement in convergence only with
relatively large amounts of memory devoted to the {\tt restart} parameter, i.e.,
 ${\tt restart}(A_{\rm R}AP^{-1})\ge 2r$. Moreover, 
preconditioned GMRES with ${\tt restart}(A_{\rm R}AP^{-1})=r/2$ did not converge 
for $\varepsilon/\varepsilon_{\rm b}=14$ at all (we set this point arbitrarilly high in 
the plot). 
On the other hand, the maximum speed-up for ${\tt restart}=2r$
was $4.1$ times (ninimum $1.8$ times), while for ${\tt restart}=10r$ the maximum 
speed-up was $11.4$ times (minimum $1.6$ times).
From the previous experiments we may expect that the speed-up will be
larger with larger objects.

\begin{figure}[t!]
\vspace*{0.5cm}
\centering
   \includegraphics[scale=1]{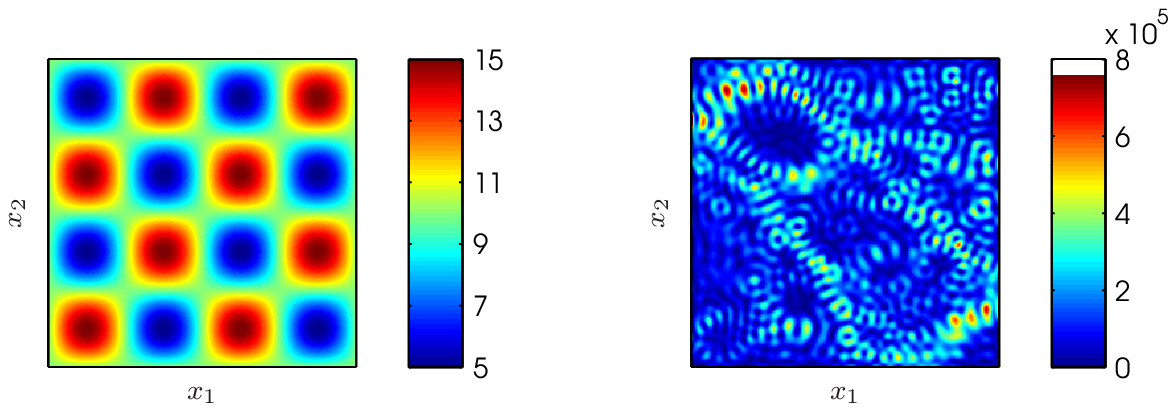}
\caption{Scattering of a TE plane wave on an inhomogeneous object with continuously varying permittivity.
Left -- permittivity profile (side length $a/\lambda_{\rm b}=4$). 
Right -- computed electric field intensity $\vert E_1\vert^2+\vert E_2\vert^2$.
Preconditioning gives $14$ times speed-up for this problem.}
\label{scatterer_fields}
\end{figure}

As an example of a systematic application of the present preconditioner, we consider
the problem, for which the DIE method is most suited, i.e., scattering on an object with continuously
varying permittivity. The steps we suggest to follow are:
\begin{enumerate}
 \item{Determine the total available memory: $M$~bytes}
 \item{Fix the discretization rule at $k$ points per smallest medium wavelength.}
 \item{Estimate $\beta$ from low-frequency matrix spectra.}
 \item{Determine the maximal number of grid points $N_{\rm max}$:}
 \begin{itemize}
           \item{If $a$ -- size of the computational domain $D$, and $N$ -- number of grid points, then
           $N=\alpha^2(a/\lambda_{\rm b})^2+2\alpha(a/\lambda_{\rm b})+1$, where\\
           $\alpha=k\max\{\sqrt{(\varepsilon'(\bx)/\varepsilon_{\rm b})}\}$, $\bx\in D$.}
           \item{Memory needed for the system, right-hand side, and the unknown vector: 
           $M_{A}=720 N$~bytes (one complex number -- 16 bytes).}
           \item{Memory required for deflation and GMRES:\\
           $M_{\rm Prec}=16r^2+48Nr+48N{\tt restart}$~bytes, where
           $r=3N\beta$ and \\ ${\tt restart}=xr$,
           i.e., $M_{\rm Prec}=144(\beta^2+\beta+x\beta)N^2$~bytes.}
           \item{The largest affordable grid follows from $M_{A}+M_{\rm Prec}\approx M$:\\
          $N_{\rm max}\approx[\sqrt{900+M(\beta^2+\beta+x\beta)}-30]/[12(\beta^2+\beta+x\beta)]$.}
 \end{itemize}
 \item{Determine the maximum affordable object size (smallest wavelength) via 
 $[a/\lambda_{\rm b}]_{\rm max}=(\sqrt{N_{\rm max}}-1)/\alpha$.}
\end{enumerate}
The scatterer is depicted in \fig{scatterer_fields}~(left). The relative permittivity function profile is given by
\begin{align}\label{eq:Permittivity_real}
\varepsilon(\bx)/\varepsilon_{\rm b}=10+5\sin(4\pi x_1/a)\sin(4\pi x_2/a),
\end{align}
where the origin of the coordinate system is in the middle of the computational domain.
We calculated $\beta$ for this scatterer to be approximately $0.0011$. Having $M=4$~GB of memory available and choosing 
$k=15$ and ${\tt restart}(A_{\rm R}AP^{-1})=8r$, from the above algorithm we determine the 
largest affordable size with this permittivity as $[a/\lambda_{\rm b}]_{\rm max}=4$. Hence,
the total number of unknowns is $3N=162867$ and the deflation parameter is $r=180$. 
The GMRES was restarted at ${\tt restart}(A_{\rm R}AP^{-1})=1440$ and ${\tt restart}(A)=1620$.
The results of this test are as follows: original system -- $162692$ seconds of CPU time, 
preconditioned system -- $11620$ seconds. The speed-up in pure `on-line' calculations
is $14$ times, as we would expect from our previous experiments with the homogeneous test object. 
The off-line work with the tolerance of {\tt eigs} set to $10^{-4}$ took only $775$ seconds,
so that even with this time included we get a $13.1$ times speed-up.
The computed electric field is shown in \fig{scatterer_fields}~(right).

The CPU timing is provided here merely for illustration and should not be used 
for direct comparison against the local methods, which was not the purpose of this paper either. 
For the particular problems considered above one could, for instance, substantially accelerate the calculations
and free some valuable memory by computing the electric field only. On the other hand, the presented CPU timing 
data give us a good relative estimate of what to expect with objects having both electric and magnetic contrasts
so that both electric and magnetic fields must indeed be computed at the same time.

Finally we should mention that, while we have been focusing here on the 
restarted GMRES, obviously, any other memory-efficient solver, which works with the original system $Au=b$, 
could also be applied to the preconditioned system $A_{\rm R}AP^{-1}v=A_{\rm R}b$. 
We have considered, for instance, Bi-CGSTAB, restarted GCR, and IDR. 
Sometimes we could achieve faster convergence by tuning the restart method of GCR, although, this did 
not happen with all types of scatterers (only with a homogeneous lossless one). Moreover, with some 
scatterers (inhomogeneous with losses) GCR and other algorithms would simply stagnate. Thus, for the moment,
the restarted GMRES remains the most robust solver for the present problem, whereas other algorithms require
more tuning and a more systematic study.

\section{Conclusions}
The analysis of the domain integral operator of the transverse electric scattering presented
here helps to understand the reasons behind the extremely slow convergence of iterative methods often 
observed with this method. Unlike the TM case, which is equivalent to the Helmholtz equation 
and can be easily accelerated by deflating the largest eigenvalues,
the TE case showed no improvement with this type of preconditioner. Direct numerical comparison of the 
TE and TM matrix spectra indicates that the nontrivial essential spectrum of the TE operator results in 
extended and very dense eigenvalue clusters, that cause stagnation of the deflation process with 
high-contrast objects. 
We have derived an analytical multiplicative regularizer  
and demonstrated on various test scatterers that its discretized version robustly reduces the extent 
of the troublesome dense clusters. One can view the obtained regularizer as a generalization of 
the Calder{\'o}n identity, previously applied in the preconditioning of boundary integral equations. 
Notably, the discretized regularizer does not require any extra computer memory and its action on a 
vector can be computed at FFT speed.
Applying deflation of the largest-magnitude eigenvalues on the regularized system we could achieve up to $46.6$ times 
acceleration of the restarted GMRES with relatively large memory, and up to $5.1$ times
with modest memory. The off-line work typically takes
only a small fraction of the total time, since one can apply a very rough tolerance in the
{\tt eigs} algorithm when recovering the largest-magnitude eigenvalues/eigenvectors.
Somewhat surprisingly, this preconditioner showed a tendency 
to become more efficient (on average) with the increase in the object size and permittivity.

\section*{Acknowledgments}
The authors are greatful to Prof. Kees Vuik (Delft University of Technology) for 
his comments and stimulating discussions.
 
\bibliographystyle{siam}

\bibliography{cleanbib}
%

\end{document}